\documentstyle[a4,amsfonts,12pt]{article}

\newcommand{\mb}{\mbox}
\newcommand{\beq}{\begin{equation}}
\newcommand{\eeq}{\end{equation}}
\newcommand{\ueberschrift}{\bigskip\goodbreak\noindent\bigskip}
\newcounter{theabsatz}
\newcommand{\absatz}[1]{\stepcounter{theabsatz} \ueberschrift
               {\large \bf \arabic{theabsatz}. {#1}} \setcounter{equation}{0}}

 \newtheorem{theor}{Theorem}
 
 \newtheorem{lem}{Lemma}
 
\newcommand{\z}{\zeta }
\begin{document}
\newcommand{\ext}{\mbox{ext\,}}
\newcommand{\diam}{\mbox{diam\,}}

\parindent 10 pt
\parskip 8pt plus 4pt
\jot 10pt

\abovedisplayskip 8pt plus 1pt \belowdisplayskip 8pt plus 1pt

\newcommand{\C}{{\mathbb{C}}}
\newcommand{\rbox}{$\:\:$ \raisebox{-1ex}{$\:\Box\:$}}
\newcommand{\CC}{{\overline{\mathbb{C}}}}
\newcommand{\DD}{{\overline{\mathbb{D}}}}
\newcommand{\D}{{\mathbb{D}}}
\newcommand{\R}{{\mathbb{R}}}
\newcommand{\He}{{\mathbb{H}}}
\newcommand{\T}{{\mathbb{T}}}
\newcommand{\cP}{{\cal P}}
\newcommand{\N}{{\mathbb{N}}}
\newcommand{\PP}{{\mathbb{P}}}
\newcommand{\p}{\preceq}
\newcommand{\s}{\succeq}
\newcommand{\En}{{\mathcal E}_n}
\newcommand{\ov}{\overline}
\newcommand{\de}{\delta}
\newcommand{\Ga}{\Gamma}
\newcommand{\ga}{\gamma}
\newcommand{\la}{\lambda}
\newcommand{\kap}{\mb{cap\,}}
\newcommand{\be}{\beta}
\newcommand{\Om}{\Omega}
\newcommand{\om}{\omega}
\newcommand{\al}{\alpha}
\newcommand{\ve}{\varepsilon}

\begin{center}
{\large \bf     Bernstein Polynomial Inequality 
on a Compact Subset of the Real Line}
\\[3ex] {\bf Vladimir Andrievskii}\\[3ex]

{\it  Department of Mathematical Sciences, Kent State University,
 Kent, OH 44242}\\[4ex]

 \end{center}

 \vspace*{2cm}

Running head: Bernstein inequality

 \vspace*{2cm}

Mailing address:

\noindent V.V. Andrievskii\\   Department of Mathematical Sciences\\
Kent
State University \\ Kent, OH 44242, USA\\[4ex]

{\it E-mail address}:\\ andriyev@math.kent.edu\\

{\it Phone}: (330) 672 9029

\newpage

\begin{center}
 {\bf Abstract}
 \end{center}
  \vspace*{0.5cm}
 We prove an analogue of the classical Bernstein polynomial inequality
on a compact subset $E$ of the real line.
The Lipschitz continuity of the Green function for the complement of
$E$ with respect to the extended complex plane and the differentiability
at a point of $E$ of a special, associated with $E$, conformal mapping
of the upper half-plane onto the comb domain play crucial role in our 
investigation.

 \vspace*{2cm}

 {\it Key Words:} Polynomial inequality, Bernstein inequality, Green's function.

 {\it AMS classification:} 30C85,  31A15, 41A17.

\newpage

\absatz{Introduction and the main result}

Let $E\subset\R$ be a non-polar compact set, i.e., there exists
the Green function $g_\Om(z)=g_\Om(z,\infty)$ of $\Om:=\CC\setminus E$
with pole at infinity. Denote by $\PP_n,n=1,2,\ldots$ the set of all
(real) polynomials of degree at most $n$ and let $||\cdot||_E$ be the supremum
norm on $E$. The classical Bernstein inequality states that for $p_n\in \PP_n$,
\beq\label{1.1}
|p_n'(x)|\le\frac{n}{\sqrt{1-x^2}}||p_n||_{[-1,1]},\quad
x\in(-1,1).
\eeq
Recently, Totik \cite{tot18} found the conditions on $E$ and $x_0\in E$
under which the analogue of (\ref{1.1}), i.e., the inequality
\beq\label{1.2}
|p_n'(x_0)|\le c(x_0,E)n||p_n||_E
\eeq
is true.
In particular, it follows directly from \cite[Theorem 3.3]{tot01} and
\cite[Lemma 3]{and18} (for details, see \cite[Theorem 2]{tot18})
that (\ref{1.2}) is equivalent to the fact that $g_\Om$ is {\it Lipschitz continuous}
at $x_0$, i.e.,
\beq\label{1.3}
\limsup_{\Om\ni z\to x_0}\frac{g_\Om(z)}{|z-x_0|}<\infty.
\eeq
For the properties of $E$ with  (\ref{1.3}) we refer
the reader to \cite{cartot, tot06, cargar, and05, and08} and the
 many references therein. 

Baran \cite[p. 489]{bar} and Totik \cite[Theorems 3.2 and 3.3]{tot01} independently found
the exact value of the constant
$c(x_0,E)$ in the case where $x_0$ is the interior (with respect to $\R$) point of $E$.
 Our objective is to extend their result to the general case
of $E$ and $x_0$ satisfying (\ref{1.3}). 
To achive this, we prove Theorem \ref{th1}, 
 which is of  independent interest.
\begin{theor} \label{th1}
Let $E$ and $x_0\in E$ satisfy (\ref{1.3}), then

(i) there exists a finite nonzero normal derivative
$$
\frac{\partial g_\Om(x_0)}{\partial n}=:h(x_0,E).
$$
Also if  $E$ is regular for the Dirichlet problem in $\Om$, then

(ii)
for 
$E_\de:=E\cup[x_0-\de,x_0+\de],\de>0$,
\beq\label{1.31}
\lim_{\de\to 0^+}h(x_0, E_\de)=h(x_0,E).
\eeq
\end{theor}
A straightforward application of Theorem \ref{th1} and \cite[Theorems 3.2 and 3.3]{tot01} 
yields the following statement.
\begin{theor} \label{th2}
Let $E$ and $x_0\in E$ satisfy (\ref{1.3}), then

(i) for $p_n\in\PP_n$,
$$
|p_n'(x_0)|\le h(x_0,E) n||p_n||_E.
$$
Also if $E$ is regular for the Dirichlet problem in $\Om$,
then 

(ii) for any
$\ve>0$ there exists $N=N(\ve)$ such that for $n>N$,
there is $p_n\in \PP_n,p_n\not\equiv 0$ satisfying
$$
|p_n'(x_0)|\ge(1-\ve)h(x_0,E) n||p_n||_E.
$$
\end{theor}
{\bf Proof}.
Since by Theorem \ref{th1}(i), \cite[Theorem 3.2]{tot01}, and by the monotonicity of the Green function with 
respect to the region for $p_n\in\PP_n$,
$$
|p_n'(x_0|\le h(x_0,E_\de)n||p_n||_{E_\de}\le
h(x_0,E)n||p_n||_{E_\de},
$$
taking the limit as $\de\to 0$, we obtain (i).

Next, according to Theorem \ref{th1}(ii), for any $\ve>0$ there is
$\de=\de(\ve)>0$ with
$$
h(x_0,E_\de)\ge \left(1-\frac{\ve}{2}\right)h(x_0,E).
$$
Moreover, by \cite[Theorem 3.3]{tot01}, there exists  $N=N(\de,\ve)$
such that for any $n>N$, there is $p_n\in\PP$  satisfying
$$
|p_n'(x_0)|\ge\left(1-\frac{\ve}{2}\right)h(x_0,E_\de)||p_n||_{E_\de}
\ge(1-\ve)h(x_0,E)||p_n||_{E},
$$
which yields (ii).

\hfill\rbox

It was shown in \cite[Corollary 2.3]{tot14}
that the Bernstein factor $h(x_0,E)n$ found in Theorem \ref{th2}
can also be stated in another form as a limit
of equilibrium densities.

\absatz{Proof of Theorem \ref{th1}(i)}

Let $E$ and $x_0\in E$ satisfy (\ref{1.3}).
 Consider the analytic in $\He:=\{z:\Im z>0\}$ function
\beq\label{2.1}
f(z)=f(z,E):= i \left(\int_E\log(z-\z)d\mu_E(\z)-\log\kap E\right),
\quad z\in\He,
\eeq
where $\kap E$ is the logarithmic capacity of $E$ and $\mu_E$ is
the equilibrium measure for $E$ (see  \cite{ran} or \cite{saftot}
for the basic notions of the potential theory).

Note that
\beq\label{2.pu2}
\Im f(z)=g_\Om(z)>0,\quad z\in\He.
\eeq
It is well known and easy to see
that $S:=f(\He)\subset\He$ is a ``comb domain", i.e,
the boundary of $S$ consists of $\R$ and at most a countable number of
closed vertical intervals with one endpoint on $\R$.
Conformal mappings of $\He$ onto such domains
play significant role in several areas of analysis (see, for example \cite{ereyud}).
We need only some elementary properties of $S$ which are discussed below.

Repeating the reasoning from \cite[pp. 222-223]{and04} one can show that $f$ is a univalent function with 
the following properties.

First, since by (\ref{1.3}) $x_0$ is a regular point,
according to  the Monotone Convergence Theorem (see \cite[p. 21]{rud})
there exists
$$
\lim_{y\to 0^+}f(x_0+iy)= -\pi\mu_E(E\cap (-\infty,x_0])=:w_0.
$$
Second, for any $z\in S,$
$$
\{\z\in S:\Re\z=\Re z,\Im\z>\Im z\}\subset S.
$$
Third, function $f$ satisfies the boundary correspondence $f(\infty)=\infty$.

Moreover,
(\ref{1.3}) and (\ref{2.1}) imply that
\beq\label{2.23}
\{w=w_0+i\eta:\eta>0\}\subset S.
\eeq
Indeed, assume that (\ref{2.23}) does not hold, then for some $d>0$
we have \newline
$[w_0,w_0+id]\subset\ov{\He}\setminus S$. To get a
contradiction,
 we use the notion of the  module
of a family of curves.
We refer to
  \cite[Chapter 4]{ahl}, \cite[Chapter IV]{garmar} or
\cite[pp. 341-360]{andbla} for the definition and basic properties of
the module such as conformal invariance, comparison principle,
composition laws, etc. We use these properties without further 
references.

Consider the Jordan curve 
$\ga:=f((x_0,x_0+i])\cup\{w_0\}$.
Without loss of generality we can assume that $\ga$ ``approaches $w_0$ 
from the right", i.e., there
is $0<\ve<d$ such that
$$
\ga\cap \{w:|w-w_0|=\ve,\pi/2<\arg(w-w_0)<\pi\}=\emptyset.
$$
Denote by $D\subset\{w\in\He:0<\arg(w-w_0)<\pi/2\}$
a Jordan domain bounded by $[w_0,w_0+id]$, a subarc
of $\ga$ and a subarc of a circle $\{w:|w-w_0|=d\}$.
For $0<\de<d$, denote by $w_\de\in D$
any point satisfying
$$
|w_\de-w_0|=\de,\quad 
\frac{\pi}{4}<\arg(w_\de-w_0)<\frac{\pi}{2},
$$
$$
\left\{w:|w-w_0|=\de,\arg(w_\de-w_0)<\arg(w-w_0)<\frac{\pi}{2}\right\}\subset D.
$$
Let $z_\de:=f^{-1}(w_\de)$. We also assume that $\de$ is so small that
$|z_\de-x_0|<1$.
Let $z_1:=x_0+iu, w_1:=f(z_1)$, where a sufficiently large but fixed number $u\ge1$ is chosen as follows. 
Denote by $\Ga_1=\Ga_1(z_\de)$ the family of all crosscuts of $\He$
which separate $z_\de$ and $x_0$ from $z_1$ and $\infty$.
Let $\Ga_2=\Ga_2(w_\de)$ be the family of all crosscuts of the
quadrilateral
$$
\{w:|w_\de-w_0|<|w-w_0|<d, 0<\arg(w-w_0)<\pi/2\}
$$ 
which separate its boundary circular components.
We choose $u$ so large that for each $\ga_2\in\Ga_2$ there exists
$\ga_1\in f(\Ga_1)$ with the property  $\ga_1\subset\ga_2$ ,
which yields that $m(\Ga_2)\le m(f(\Ga_1)).$

According to \cite[(2.5)]{and18}, we obtain
$$
\frac{1}{\pi}\log\frac{u}{|z_\de-x_0|}+2
\ge m(\Ga_1)= m(f(\Ga_1))\ge m(\Ga_2)
=\frac{2}{\pi}\log\frac{d}{\de},
$$
which implies
$$
g_\Om(z_\de)=\Im w_\de>\frac{\de}{2}\ge
\frac{d}{2e^\pi\sqrt{u}}|z_\de-x_0|^{1/2}.
$$
Therefore,
$$
\frac{g_\Om(z_\de)}{|z_\de-x_0|}\to \infty\quad \mb{as }\de\to 0,
$$
which contradicts (\ref{1.3}). This completes the proof
of (\ref{2.23}).

Next, for $l>0$ and $v>0$, consider 
the quadrilateral
$$
Q_{l,v}:=\{ z\in\He:v<|z|<v\sqrt{1+l^2}, \Re z<v\}.
$$
Denote by $\Ga_{l,v}$ the family of all crosscuts of $Q_{l,v}$ that separate its 
circular boundary
components. 

We claim that there exists a positive constant $c=c(l)$ such that
\beq
\label{2.0}
m(\Ga_{l,v})-\frac{1}{\pi}\log\sqrt{1+l^2}\ge c.
\eeq

Indeed, applying the  transformation $z\to z/v$ and using conformal
invariance of the module, we can reduce the proof of (\ref{2.0}) to
the case where
 $v=1$.
Since for $r\in(1,\sqrt{1+l^2})$, 
$$
\{\theta:re^{i\theta}\in Q_{l,1}\}=(\cos^{-1}(1/r),\pi),
$$
for the module  of $\Ga_{l,1}$ we have
\begin{eqnarray*}
&& m(\Ga_{l,1})-\frac{1}{\pi}\log\sqrt{1+l^2}\\
&\ge&
\int_1^{\sqrt{1+l^2}}\frac{dr}{r(\pi-\cos^{-1}(1/r))}
-\frac{1}{\pi}\log\sqrt{1+l^2}\\
&=&
\int_1^{\sqrt{1+l^2}}\frac{\cos^{-1}(1/r)dr}{\pi r(\pi-\cos^{-1}(1/r))}\\
&\ge& 
\frac{1}{\pi^2\sqrt{1+l^2}}
\int_1^{\sqrt{1+l^2}}\cos^{-1}(1/r)dr=:c,
\end{eqnarray*}
which yields (\ref{2.0}).

In the proofs of some of the lemmas below we use the following family of curves
and its module.
 For $w_1=w_0+i\eta$ and $w_2=w_0+i$ with $0<\eta<1$ we introduce the family
$\Ga_0=\Ga_0(\eta)$
of all crosscuts of $S$ that separate  points $w_0$ and
$w_1$ from $w_2$ and $\infty$ in $S$. 
Let $z_k:=f^{-1}(w_k), k=1,2$.
According to \cite[(2.5)]{and18}, if $|z_1-x_0|<|z_2-x_0|,$ then for the module of $\Ga_0$ we have
\beq\label{2.6}
m(\Ga_0)=m(f^{-1}(\Ga_0))\le\frac{1}{\pi}\log\left|\frac{z_2-x_0}{z_1-x_0}\right|
+2.
\eeq
Following \cite[p. 173]{garmar}, for $\be\in(0,\pi/2)$ and $\ve>0$, define the
{\it truncated cone}
$$
S_\be^\ve=S_\be^\ve(w_0):=\{w:|\arg(w-w_0)-\pi/2|<\ve,
0<|w-w_0|<\ve\}.
$$
We say that $\partial S$
 has an {\it inner tangent} (with inner normal $i$) at $w_0$ if for every 
 $\be\in(0,\pi/2)$ there is an $\ve=\ve(\be)>0$ so that
$S_\be^\ve\subset S$.
\begin{lem}\label{lem1}
$\partial S$ has an inner tangent at $w_0$.
\end{lem}
{\bf Proof}.
Assume that $\partial S$ does not have an inner tangent at $w_0$.
Then, there exist $\ve>0$ and a sequence of real numbers $\{b_n\}_1^\infty$ such that
$$
d_n:=|b_n-w_0|<1,\quad \lim_{n\to \infty}b_n=w_0,
$$
\beq\label{2.4}
[b_n,b_n+id_n]\subset\ov{\He}\setminus S,
\eeq
\beq\label{2.5}
\sqrt{1+\ve^2}d_1<1,\quad
\sqrt{1+\ve^2}d_{n+1}<d_n.
\eeq
For $r>0$, denote by $\ga(r)=\ga(w_0,S,r)\subset\{w:|w-w_0|=r\}$
the crosscut of $S$ which has a nonempty intersection with the ray
$\{w\in \He: \Re w=w_0\}$. For $0<r<R$, denote by $D(r,R)=D(w_0,S,r,R)\subset
S$ the bounded simply connected domain whose boundary consists of $\ga(r),\ga(R)$,
and two connected parts of $\partial S$. Let $m(r,R)=m(w_0,S,r,R)$ be the module of the family 
$\Ga(r,R)=\Ga(w_0,S,r,R)$ of all
crosscuts of $D(r,R)$ which separate  circular arcs $\ga(r)$ and $\ga(R)$ in $D(r,R)$.
Note that
\beq\label{2.1n}
m(r,R)\ge\frac{1}{\pi}\log\frac{R}{r},\quad 0<r<R.
\eeq
Moreover, according to (\ref{2.0}), (\ref{2.4}), and (\ref{2.5}),
$$
m(d_n, d_n\sqrt{1+\ve^2})-\frac{1}{\pi}\log\sqrt{1+\ve^2}\ge c=c(\ve)>0.
$$
Therefore, by virtue of (\ref{2.1n}) we obtain
\begin{eqnarray*}
&&m(\Ga_0(d_n))\ge m(d_n,1)\\
&\ge&
m(\sqrt{1+\ve^2}d_1,1)+\sum_{k=1}^n m(d_k,\sqrt{1+\ve^2}d_k)+
\sum_{k=2}^n m(\sqrt{1+\ve^2}d_k,d_{k-1})\\
&\ge&
\frac{1}{\pi}\log\frac{1}{d_n}+
\sum_{k=1}^n \left( m(d_k,\sqrt{1+\ve^2}d_k)-\frac{1}{\pi}\log\sqrt{1+\ve^2}
\right)\\
&\ge&
\frac{1}{\pi}\log\frac{1}{d_n}+nc.
\end{eqnarray*}
Comparing the last inequality with (\ref{2.6}) for $\eta=d_n$,
 we have
$$
\frac{g_\Om(z_1)}{|z_1-x_0|}\ge\frac{e^{\pi(cn-2)}}{|z_2-x_0|}\to 
\infty\quad \mb{as }n\to\infty,
$$
which contradicts (\ref{1.3}).

\hfill\rbox

According to Lemma \ref{lem1} and the Ostrowski Theorem
(see \cite[p. 177, Theorem 5.5]{garmar}) for every
$\be\in(0,\pi/2)$ there exist
nontangential limits
$$
\lim_{S_\be^\ve(x_0)\ni z\to x_0} f(z)=w_0,
$$
$$
\lim_{S_\be^\ve(x_0)\ni z\to x_0} \arg\frac{f(z)-w_0}{z-x_0}=0.
$$
\begin{lem}\label{lem2}
The function $f$ has a positive angular derivative at $x_0$, i.e.,
\beq\label{2.9}
f'(x_0):=
\lim_{S_\be^\ve(x_0)\ni z\to x_0} \frac{f(z)-w_0}{z-x_0}>0
\eeq
exists for every $\be\in(0,\pi/2)$.
\end{lem}
{\bf Proof}. Assume that (\ref{2.9}) is not true.
By the Jenkins-Oikawa-Rodin-Warschawski Theorem (see
\cite[p. 180, Theorem 5.7]{garmar}) there exists $\ve>0$ such that 
for every $\de>0$ there are  $0<s<r<\de$  satisfying
$$
m(s,r)-\frac{1}{\pi}\log\frac{r}{s}\ge\ve.
$$
Therefore, we can find sequences of real numbers
$\{s_n\}_1^\infty$ and $\{r_n\}_1^\infty$
such that $0<r_n<s_{n-1}<r_{n-1}<1$ and
$$
\lim_{n\to\infty}s_n=\lim_{n\to\infty}r_n=0,
$$
$$
m(s_n,r_n)-\frac{1}{\pi}\log\frac{r_n}{s_n}>\ve.
$$
Furthermore, by (\ref{2.1n})
\begin{eqnarray*}
m(\Ga_0(s_n))&\ge&m(s_n,1)\ge m(r_1,1)+
\sum_{k=1}^nm(s_k,r_k)+\sum_{k=1}^{n-1}m(r_{k+1},s_k)\\
&\ge&
\frac{1}{\pi}\log\frac{1}{s_n}+
\sum_{k=1}^n\left(m(s_k,r_k)-\frac{1}{\pi}\log\frac{r_k}{s_k}
\right)\\
&\ge&
\frac{1}{\pi}\log\frac{1}{s_n}+n\ve.
\end{eqnarray*}
Comparing the last inequality with (\ref{2.6}) for $\eta=s_n$, we obtain
$$
\frac{g_\Om(z_1)}{|z_1-x_0|}\ge\frac{1}{|z_2-x_0|}e^{\pi(\ve n-2)}\to
\infty\quad \mb{as }n\to\infty,
$$
which contradicts (\ref{1.3}).

\hfill\rbox

Note that Lemma \ref{lem2} and (\ref{2.pu2}) imply Theorem \ref{th1}(i).

\absatz{Proof of Theorem \ref{th1}(ii)}

In this section we assume that $E$ is a regular set satisfying (\ref{1.3}).
Hence, the extension of the Green function by letting
$g_\Om(z):=0,z\in E$ produces a continuous function in $\C$.
\begin{lem}\label{lem3}
For the Green function we have
\beq\label{3.1}
\int_{-\infty}^\infty\frac{g_\Om(x) dx}{(x-x_0)^2}<\infty.
\eeq
\end{lem}
{\bf Proof}.
Consider $E^*:=\{1/(x-x_0):x\in E\}\subset \R\cup\{\infty\}\setminus
\{0\}.$ Note that
$$
g_\Om(z)=g_{\CC\setminus E^*}(1/(z-x_0),0),\quad z\in\CC.
$$
Denote by $\cP_\infty$ the cone of positive harmonic functions on $\C\setminus
E^*$ which have vanishing boundary values at every point of $E^*\setminus\{\infty\}.$
Since (\ref{1.3}) implies
$$
\limsup_{\C\setminus E^*\ni \z\to\infty}g_{\C\setminus E^*}(\z,0)|\z|<\infty,
$$
according to \cite[Theorem 1]{and08} or \cite[Theorem 1]{cargar}
the dimension of $\cP_\infty$ is $2$. Therefore, by virtue of
\cite[Theorem 3]{kar} (cf. \cite[Section VIII A.2]{koo}), 
there exists a Phragm\'en-Lindel\"of function for $\C\setminus E^*$
which yields 
$$
\int_{-\infty}^\infty
g_{\C\setminus E^*}(t,0)dt<\infty
$$
and (\ref{3.1}) follows.

\hfill\rbox

According to Lemma \ref{lem3}
and Monotone Convergence Theorem (see\cite[p. 21]{rud}),
for any $\ve>0$ there exists $\de_1=\de_1(\ve)>0$
with the property
\beq\label{3.2}
\int_{x_0-\de_1}^{x_0+\de_1}\frac{g_\Om(x) dx}{(x-x_0)^2}<\ve.
\eeq
Next, let $0<\de=\de(\ve,\de_1)<\de_1$ be such that for $x_0-\de\le x\le x_0+\de$
the inequality
$g_\Om(x)<\ve\de_1$ holds. Consider the function
$$
u_\de(z):= g_\Om(z)-g_{\Om_\de}(z),\quad z\in \Om_\de:= \CC\setminus E_\de,
$$ 
which is harmonic in $\Om_\de$.

By the maximum principle
$$
u_\de(x)\le\ve\de_1,\quad x\in \R.
$$
Moreover, (\ref{3.2}) implies that
$$
\int_{x_0-\de_1}^{x_0+\de_1}\frac{u_\de (x) dx}{(x-x_0)^2}\le
\int_{x_0-\de_1}^{x_0+\de_1}\frac{g_\Om(x) dx}{(x-x_0)^2}<\ve.
$$
Applying the Poisson formula (see \cite[p. 4]{garmar}), for $y>0$ we
obtain
\begin{eqnarray*}
u_\de(x_0+iy)&=&
\frac{y}{\pi} \int_{-\infty}^\infty\frac{u_\de(x)dx}{(x-x_0)^2+y^2}\\
&\le&
\frac{y}{\pi} \left(\int_{|x-x_0|\ge\de_1}\frac{u_\de(x)dx}{(x-x_0)^2}
+ \int_{|x-x_0|<\de_1}\frac{g_{\Om_\de}(x)dx}{(x-x_0)^2}\right)\\
&\le&\frac{y}{\pi}\left(\frac{2}{\de_1}\ve\de_1+\ve\right)<y\ve.
\end{eqnarray*}
Therefore,
\beq\label{3.3}
\frac{\partial u_\de(x_0)}{\partial n}:=\lim_{y\to 0^+}
\frac{u_\de(x_0+iy)}{y}\le\ve.
\eeq
Monotonicity of the Green function and (\ref{3.3}) yield
$$
\frac{\partial g_{\Om_\de}(x_0)}{\partial n}\le
\frac{\partial g_{\Om}(x_0)}{\partial n}\le
\frac{\partial g_{\Om_\de}(x_0)}{\partial n}+\ve,
$$
which completes the proof of (\ref{1.31}).

\absatz{Acknowledgements}
The
author is grateful to  M. Nesterenko, V. Totik and A. Volberg
  for their helpful comments.


\begin{thebibliography}{99}


\bibitem{ahl}
L.V. Ahlfors, Conformal Invariants, McGraw-Hill, New York, 1973.
 


\bibitem{and04}
 V.V. Andrievskii,   The highest smoothness of the
Green function implies the highest density of a set,  Ark.
Mat.
 42  (2004)   217--238.


\bibitem{and05}
V.V. Andrievskii, On sparse sets with the Green function
of the highest smoothness, Comp. Met. and Fun. Theory, 
5 (2005) 301--322.

\bibitem{and08}
V.V. Andrievskii, Positive harmonic functions on Denjoy domains in
the complex plane,  J. d'Analyse Math.  104 (2008) 83--124.

\bibitem{and18}
V.V. Andrievskii, Polynomial approximation on a compact subset of the
real line, Journal of  Approximation Theory 230 (2018) 24--31.


\bibitem{andbla} V.V. Andrievskii, H.-P. Blatt, Discrepancy of Signed Measures and
      Polynomial Approximation,
       Springer-Verlag, Berlin/New York, 2002.
			
			
\bibitem{bar}
M. Baran. Complex equilibrium measure and Bernstein type theorems for
compact sets in $\R^n$, Proc. Amer. Math. Soc. 123(2) (1995) 485--494.			

  \bibitem{cartot}
 L. Carleson, V. Totik,  H\"older continuity of Green's
 functions,  Acta Sci. Math. (Szeged)  70  (2004)   557--608.

 \bibitem{cargar}
 T. Carroll, S. Gardiner,   Lipschitz continuity of the
 Green function in Denjoy domains, Ark. Mat. 46 (2008)
271--283.

\bibitem{ereyud}
A. Eremenko, P. Yuditskii, Comb functions,  Contemp. Math. 578 (2012) 99--118. 

\bibitem{garmar}
J.B. Garnet, D.E. Marshal, Harmonic Measure, 
Cambridge University Press, Cambridge, 2005.

\bibitem{kar}
P.P. Kargaev, Existence of the Phragm\'en-Lindel\"off
function and certain quasianalyticity conditions,
Zap. Nauchn. Sem. Leningrad. Otdel. Mat. Inst. Steklov.
(LOMI) 126 (1983) 97--108
(Russian); English translation in J. Soviet Math. 27 (1984)
2486-2495.

\bibitem{koo}
 P. Koosis, The Logarithmic Integral, I, Cambridge University Press,
Cambridge, 1988.




\bibitem{ran}
 T. Ransford,  Potential Theory in the Complex Plane,
   Cambridge University Press, Cambridge, 1995.

\bibitem{rud}
W. Rudin, Real and Complex Analysis, 3rd ed., McGraw-Hill,
New York, 1987.

\bibitem{saftot}
 E.B. Saff ,  V. Totik,    Logarithmic Potentials with External
Fields, Springer-Verlag,  New York/Berlin, 1997.

\bibitem{tot01}
 V. Totik, Polynomial inverse images and polynomial inequalities, Acta
Math. 187 (2001) 139--160.



\bibitem{tot06}
 V. Totik,   Metric Properties of Harmonic Measures,
Mem. Amer. Math. Soc. 184  (2006), no 867.

\bibitem{tot14}
V. Totik, Bernstein- and Markov-type inequalities for trigonometric
polynomials on general sets, International Mathematics Research Notices, rnu030,
35 pages.
doi:10.1093/imrn/rnu030

\bibitem{tot18}
V. Totik, Reflection on a theorem of V. Andrievskii,
manuscript, 2018.


\end{thebibliography}
\end{document}